\documentclass[12pt,a4]{article}
\usepackage{amssymb,amscd,amsmath,amstext,amsfonts}
\usepackage{amsthm}
\usepackage{array}
\usepackage{xspace,rotating,color}
\usepackage{fancyhdr}
\usepackage{makeidx}

\newtheorem{thm}{Theorem}[section]
\newtheorem{prop}[thm]{Proposition}
\newtheorem{defn}[thm]{Definition}

\newtheorem{rem}[thm]{Remark}

\newtheorem{ex}[thm]{Example}

\def\BB{\mathcal B}
\def\CC{\mathbb C}
\def\DD{\mathbb D}
\def\FF{\mathcal F}
\def\GF{\FF}
\def\LL{\mathcal L}
\def\MM{\mathcal M}
\def\RR{\mathbb R}
\def\TT{\mathbb T}
\def\BMOA{\rm BMOA}
\def\fix{\mathop{\rm fix}\nolimits}
\def\Hol{\mathop{\rm Hol}\nolimits}
\def\Im{\mathop{\rm Im}\nolimits}
\def\re{\mathop{\rm Re}\nolimits}

\def\Id{{\rm Id}}
\def\Gen{\mathop{\rm Gen}\nolimits}
\def\DW{\mathop{\rm DW}\nolimits}

\def\beginpf{\begin{proof}}
\def\endpf{\end{proof}}
\def\beq{\begin{equation}}
\def\eeq{\end{equation}}

\newcommand\W{W_{w,\phi}}

\makeindex

\begin{document}  

\title{Semigroups of weighted composition operators on spaces of holomorphic functions}

\author{I. Chalendar and J.R. Partington}

\maketitle

{Isabelle Chalendar,  Universit\'e Gustave Eiffel, LAMA, (UMR 8050), UPEM, UPEC, CNRS, F-77454, Marne-la-Vallée, France.}\\
{\tt isabelle.chalendar@univ-eiffel.fr}\\

{Jonathan R. Partington, School of Mathematics, University of Leeds, Leeds LS2 9JT, UK.}
{\tt J.R.Partington@leeds.ac.uk}\\


\section{Introduction}
This paper is based on three hours of lectures given by the first author  in the    
``Focus Program on Analytic Function Spaces and their Applications"
July 1 -- December 31, 2021,   organized by  the Fields Institute for 
Research in Mathematical Sciences.

The  goal of this paper is to give an introduction to the properties of discrete and continuous $C_0$-semigroups of (weighted) composition operators on various spaces of analytic functions. 

To that aim we detail the structure of semiflows of analytic functions on the open unit disc $\DD$ and their  generators, which  provides  information on the generators of continuous semigroups of composition operators on Banach spaces $X$ embedding in $\Hol(\DD)$, the Fr\'echet space of holomorphic functions 
on $\DD$.     

An initial  motivation for studying such semigroups is a better understanding of an specific  universal operator. Moreover, adding a weight to a composition operator is motivated by the fact that such operators describe isometries on non-Hilbertian Hardy spaces and they appear automatically when the Banach 
spaces $X$ are replaced by Banach spaces of holomorphic functions on another domain such as the right 
half-plane. 

Thanks to the analysis of spectral properties,  we deduce the asymptotic behaviour of discrete semigroups of composition operators on various Banach spaces such as the Hardy spaces, the weighted Bergman spaces, Bloch type spaces or standard weighted Bergman space of infinite order. 
As a byproduct we obtain characterization of the properties of isometry and similarity to isometry, still  for composition operators.  
 We then describe the limit at infinity for continuous semigroups of composition operators.
 
 Compactness (immediate and eventual) and analyticity of semigroups of composition operators are then considered on the Hardy space $H^2(\DD)$, even though other classes of Banach spaces may also be considered. References are given for more complete information. 
 
 Finally, we provide some perspectives  for semigroups of composition operators on $H^2(\CC_+)$, where $\CC_+$ is the right-half-plane, as well as an analysis of  semigroups of composition operators on the Fock space. The latter case can be treated in a complete way since the non-trivial semiflows involved are necessarily expressed using polynomials of degree one.

\tableofcontents
\section{Background}
\subsection{Strongly continuous semigroups of operators: definition and characterization}
We recall some of the standard facts about one-parameter semigroups of operators,
which may be found in many places, such as \cite{EN06} and \cite{pazy}.

\begin{defn}
A semigroup $(T_t)_{t \ge 0}$ of operators on a Banach space $X$ is a family
of bounded operators satisfying:\\
(i) $T_0= \Id$, the identity operator, and (ii) $T_{t+s}=T_t T_s$ for all $s,t \ge 0$.\\
If, in addition, it satisfies:\\
(iii) $T_t x \to x$ as $t \to 0^{+}$ for all $x \in X$, then it is called a {\em strongly continuous} or 
$C_0$-semigroup.
\end{defn}

A {\em uniformly continuous\/} semigroup $(T_t)_{t \ge 0}$
is one satisfying $\|T_t-\Id\| \to 0$ as $t \to 0^{+}$.

Associated with this is the notion of an {\em  infinitesimal generator}, or simply {\em generator}.
We define an (in general unbounded) operator $A$ whose domain is
\[
D(A):= \{x \in X: \lim_{t \to 0^{+}} \frac{T_t x -x}{t} \quad \hbox{exists} \},
\]
and then 
\[
Ax := \lim_{t \to 0^{+}}  \frac{T_t x -x}{t} \qquad \hbox{for} \quad x \in D(A).
\]

Moreover, the generator of a $C_0$-semigroup characterizes completely a semigroup, that is  
two $C_0$-semigroups are equal if and only if  their generators are equal.

By the uniform boundedness principle, each $C_0$-semigroup is uniformly bounded on each compact interval. As a corollary, for every $C_0$-semigroup $(T_t)_{t\geq 0}$, there exists $w\in\RR$ and $M\geq 1$ such that 
   \[ \|T_t\|\leq Me^{wt}\mbox{ for all }t\geq 0.\]

Contractive $C_0$-semigroups  are the ones for which one can take $M=1$ and $w=0$, whereas quasicontractive $C_0$-semigroups are the ones for which one can take $M=1$ and $w$ is arbitrary.  \\

The domain $D(A)$ of the generator $A$ of a $C_0$-semigroup is always dense in $X$ (and moreover $(A,D(A))$ is a closed operator).  It is then natural to characterize the  linear operators $(A,D(A))$ that are the  generator of  $C_0$-semigroups.
For semigroups of contractions   on Hilbert spaces, 
Lumer and Phillips, in 1961, provided  a beautiful criterion \cite{LP61} (see also \cite[Theorem~3.15]{EN06}). 
\begin{thm}[Lumer--Phillips]
Let $(A,D(A))$ be a linear operator with dense domain on a Hilbert space $H$. The following assertions are equivalent:
\begin{enumerate}
	\item[(i)] $A$ generates a $C_0$-semigroup of contractions;
	\item[(ii)]  there exists $\lambda>0$ such that $(\lambda\Id -A)D(A)=H$ and for all $x\in D(A)$,
	\[   \re \langle Ax,x\rangle \leq 0;\]
	\item[(iii)]  for all  $\lambda>0$ we have $(\lambda\Id -A)D(A)=H$ and for all $x\in D(A)$,
	\[   \re \langle Ax,x\rangle \leq 0.\]
\end{enumerate}
\end{thm}
Since $(T_t)_{t\geq 0}$ is a $C_0$-semigroup of quasicontractions on a Hilbert space $H$ if and only if there exists $w\geq 0$ such that $(e^{-wt}T_t)_{t\geq 0}$ is a $C_0$-semigroup of contractions,  it is then easy to deduce a version of the Lumer--Phillips Theorem for $C_0$-quasicontractions.

\begin{thm}
	Let $(A,D(A))$ be a linear operator with dense domain on a Hilbert space $H$. The following assertions are equivalent:  
	\begin{enumerate}
		\item[(i)] $A$ generates a $C_0$-semigroup of quasicontractions;
		\item[(ii)]  there exists $\lambda>0$ such that $(\lambda\Id -A)D(A)=H$ and  
		\[  \sup_{x\in D(A),\|x\|\leq 1}\re \langle Ax,x\rangle <\infty ;\]
		\item[(iii)]  for all  $\lambda>0$ we have $(\lambda\Id -A)D(A)=H$ and  
		\[   \sup_{x\in D(A),\|x\|\leq 1}\re \langle Ax,x\rangle <\infty.\]
	\end{enumerate}
\end{thm}	 
For generators of $C_0$-semigroups which are not necessarily quasicontractive, another beautiful criterion involving the growth of the resolvent is due to Hille and Yosida \cite[Theorem~3.8]{EN06}. \\
   
\begin{thm}
Let $A$ be a linear operator defined on a linear subspace $D(A)$ of the Banach space $X$, $w\in\RR$ and 
$M>0$. Then $A$ generates a $C_0$-semigroup $(T_t)_{t\geq 0}$ that satisfies 
\[ \|T_t\|\leq M e^{w t}\] if and only if
\begin{itemize}
	\item[a)] $A$ is closed and $D(A)$ is dense in $X$;
	\item[(b)] every real $\lambda > w$ belongs to the resolvent set of $A$ and for such $\lambda$ and for all positive integers $n$,
	\[ \|(\lambda \Id-A)^{-n}\|\leq {\frac {M}{(\lambda -\omega )^{n}}}.\]
\end{itemize}
\end{thm}

We shall restrict ourselves to Banach spaces $X$ of functions   that are holomorphic on a domain $\Omega$ (usually
the unit disc $\DD$ but sometimes the right half-plane $\CC_+$ or the complex plane $\CC$) and satisfying the
condition that point evaluations $f \mapsto f(z)$ are continuous for all $z \in \Omega$. Assuming this,
we have that 
norm convergence of a sequence $(f_n)$ to $f$ implies  local uniform convergence (uniform convergence
on compact subsets of $\Omega$).

Recall that for suitable $\phi: \Omega \to \Omega$ holomorphic, the   composition
operator $C_\phi: X \to X$ is defined by $(C_\phi f)(z)= f(\phi(z))$, for $f \in X$ and $z \in \Omega$ (assuming that $C_\phi$
maps $X$ boundedly into $X$, an issue which will be discussed later). Likewise, for suitable $w$ holomorphic on $\Omega$,
the weighted composition operator $W_{w,\phi}$ on $X$ is defined by $(W_{w,\phi}f)(z)=w(z)f(\phi(z))$.\\

The main theme of this paper is to study $C_0$-semigroups of (weighted) composition operators.
However, we
may also look at composition semigroups from a non-operatorial point of view (for example, as in \cite{BCD20}). 

\subsection{Analytic semiflows on a domain and models for semiflows on $\DD$}
\begin{defn}
A continuous  {\em analytic semiflow\/} on a domain $\Omega$ is a family $(\phi_t)_{t \ge 0}$ of 
holomorphic mappings from $\Omega$ to itself satisfying:\\
(i) $\phi_0 (z)=z$ for all $z \in \Omega$,\\
(ii)  $\phi_{t+s}=\phi_t \circ \phi_s$ for all $s,t \ge 0$, and\\
(iii) for all $z \in \Omega$, the mapping $t \mapsto \phi_t(z)$ is continuous on $[0,\infty)$.
\end{defn}
\begin{rem}
A family $(\phi_t)_{t \ge 0}$ of 
holomorphic mappings from $\Omega$ to itself satisfying only 
(i)  and (ii) is called an \emph{algebraic semiflow}.  Such an algebraic  semiflow $(\phi_t)_{t \ge 0}$ is continuous  on $\DD$ if and only if there exists $a\in\DD$ such that $\lim_{t \to 0} \phi'_t(a)=1$ (see \cite[Thm. 8.1.16]{BCD20}).    
\end{rem}

In this situation there exists a unique holomorphic function $G:\Omega\to\CC$ such that 
\[  \frac{\partial \phi_t(z)}{\partial t}=G(\phi_t(z)),\,\,\, z\in\Omega  \mbox{ and }t\in [0,\infty).\]
This function is called the infinitesimal generator of the analytic semiflow   $(\phi_t)_{t \ge 0}$ on $\Omega$ and 
we denote by $\Gen(\Omega)$ the set of all infinitesimal generator of analytic semiflows on $\Omega$.      \\

 There are several complete characterization of $\Gen(\DD)$ in \cite[Chap. 10]{BCD20}, which is discussed in more details in Section~\ref{sec:3.2}.

Analytic semiflows on $\DD$ can be partitioned into two classes, depending on the localization of their Denjoy--Wolff point  $\alpha$, 
discussed below in Section \ref{sec:4ab}  (see \cite{Abate89}, \cite[Chap. 2]{CM95} and \cite[Chap. 8]{BCD20}). 

If $\alpha \in \DD$, 
  by conjugating by the
automorphism $b_\alpha$, where
\[
b_\alpha(z):= \frac{\alpha-z}{1-\overline\alpha z},
\]
we may suppose without loss of generality that $\alpha=0$. In this case there
is a semiflow model
\[
\phi_t(z)= h^{-1}(e^{-ct}h(z)),
\]
where $c \in \CC$ with $\re c \ge 0$, and $h: \DD \to \Omega$ is a conformal
bijection between $\DD$ and a domain $\Omega \subset \CC$, with $h(0)=0$
and $\Omega$ is spiral-like  or star-like (if $c$ is real), in the sense that
\[
e^{-ct}w \in \Omega \quad \hbox{for all} \quad w \in \Omega \quad \hbox{and} \quad t \ge 0.
\]

If $\alpha \in \TT$,  then there exists a conformal map $h$ from $\DD$ onto a domain $\Omega$ such that $\Omega +it\subset \Omega$ for all $t\geq 0$, and    there
is a semiflow model
\[
\phi_t(z)= h^{-1}(h(z)+it).
\]

\subsection{Models for analytic flows on $\DD$}
This subsection relies heavily on Subsection 8.2 in \cite{BCD20}. 
\begin{defn}
A family  $(\phi_t)_{t\geq 0}$ of analytic selfmaps of  $\DD$ is a called a continuous (algebraic) \emph{flow} if 
\begin{itemize}
	\item[a)]  $\phi_t$ is an automorphism of $\DD$ for all $t\geq 0$;
	\item[b)]  $(\phi_t)_{t\geq 0}$ is a continuous (algebraic) semiflow.
\end{itemize}
\end{defn}
If $\phi_t$ is an automorphism for all $t\geq 0$, we can introduce the notation $\phi_{-t}:=\phi_t^{-1}$ for all $t\geq 0$ and then observe that  
{\em 
 \begin{itemize}
 	\item[c)]  $\phi_{s+t}=\phi_s\circ \phi_t$ for all $s,t\in\RR$;
 	\item[d)]   for all $z \in \DD$, the mapping $t \mapsto \phi_t(z)$ is continuous on $\RR$ if $(\phi_t)_{t\geq 0}$   is  a continuous  flow.
 \end{itemize}
 }
In fact, d) is equivalent to the continuity of $t\mapsto \phi_t\in \Hol(\DD)$ on $\RR$, where $\Hol(\DD)$ is endowed with the topology of the uniform convergence on compacta of $\DD$.  

Here is a characterization of continuous flows in the set of continuous semiflows \cite[Thm. 8.2.4]{BCD20}.
\begin{thm}
	Let $(\phi_t)_{t\geq 0}$ be  a continuous (algebraic) semiflow on $\DD$. Then it is a continuous (algebraic) flow if and only if there exists $t_0>0$ such that $\phi_{t_0}$ is an automorphism.    
\end{thm}  
The following theorem \cite[Thm. 8.2.6]{BCD20} is an explicit description of all the continuous flows on $\DD$.
\begin{thm}
Let $(\phi_t)_{t\geq 0}$ be a nontrivial continuous flow on $\DD$. Then $(\phi_t)_{t\geq 0}$ has one of the following three mutually exclusive forms:
\begin{itemize}
	\item[1)] There exists $\alpha\in\DD$ and $w\in\RR\setminus\{0\}$ such that 
	\[    \phi_t(z)=\frac{(e^{-iw t} -|\alpha|^2)z +\alpha (1-e^{-iwt})  }{\overline{\alpha}(e^{-iwt} -1)z +1-|\alpha|^2e^{-iwt}},\,\,\, t\geq 0,\,\,\, z\in\DD .\] 
	Moreover it is the unique continuous flow of elliptic automorphisms for which $\phi_t(\alpha)=\alpha$ for all $t$ and $\phi'_t(\alpha)=e^{-iwt}$. 
	\item[2)] There exist $\alpha_1,\alpha_2\in\TT$, $\alpha_1\neq \alpha_2$ and $c>0$ such that 
	 	\[    \phi_t(z)=\frac{(\alpha_2 -\alpha_1e^{ct})z +\alpha_1\alpha_2 (e^{ct}-1)  }{(1-e^{ct})z +\alpha_2e^{ct}-\alpha_1},\,\,\, t\geq 0,\,\,\, z\in\DD .\] 
	 Moreover it is the unique continuous flow of hyperbolic automorphisms for which $\phi_t(\alpha_i)=\alpha_i$ for all $t$  ($i=1,2$) and   $\phi'_t(\alpha_1)=e^{-c t}$. 
	\item[3)] There exist $\alpha\in\TT$ and $w\in\RR\setminus\{0\}$ such that 
	\[ \phi_t(z)=\frac{(1-iwt)z+iw\alpha t  }{-iw\overline{\alpha} tz +1+iwt  }.   \]
	Moreover it is the unique continuous flow of parabolic automorphisms for which $\phi_t(\alpha)=\alpha$ for all $t$    and   $\phi''_t(\alpha)=2itw\overline{\alpha}$. 
\end{itemize}
\end{thm}  

\subsection{$C_0$-semigroups of composition operators}
Clearly a semiflow $(\phi_t)_{t \ge 0}$ induces a semigroup of composition operators
(on $\DD$ these are bounded, by Littlewood's subordination theorem \cite{littlewood}),
and the following condition gives a way of testing the strong continuity.

\begin{prop}\label{prop:sg}
Let  $E$ be a dense subspace of a Banach space $X$ and $(T_t)_{t \ge 0}$
a semigroup of bounded operators on $X$ such that there exists a $\delta>0$ with 
$\sup_{0 \le t \le \delta} \|T_t\|< \infty$. Then $(T_t)$ is a $C_0$-semigroup on $X$
if and only if $\lim_{t \to 0} \|T_t f -f\|_X=0$
for all $f \in E$.  
\end{prop}

In particular, if the polynomials are dense in $X$ then it is enough 
to check that $\lim_{t \to 0} \|T_t e_n - e_n \|_X = 0$ for all $n=0,1,2,\ldots$,
where $e_n(z)=z^n$.

\beginpf
Clearly, the ``only if'' condition holds. Conversely, if $\lim_{t \to 0} \|T_t f -f\|_X=0$
for all $f \in E$, 
let $M:=\sup_{0 \le t \le \delta} \|T_t\|$, and
let $\epsilon>0$ be given  and $f \in X$. 
We may find a $p \in E$ such that $\|f-p\|<\dfrac{\epsilon}{2(M+1)}$. Then
\begin{align*}
\|T_t f-f\| &\le \|T_t f-T_t p\|+ \|T_t p- p\|+ \|p-f\| \\
& \le M \|f-p\| + \epsilon/2 + \|p-f\| < \epsilon
\end{align*}
for sufficently small $t$.
\endpf

\subsection{Spaces on which semigroups are not $C_0$}

In Proposition \ref{prop:sg} we have seen a sufficient condition for a semiflow to induce a $C_0$-semigroup of composition
operators, and in the Hardy, Dirichlet and Bergman spaces we do indeed arrive at such a semigroup.

Recently, Gallardo-Guti\'errez, Siskakis and Yakubovich \cite{GSY21} have   shown that weighted composition operators
$(W_{w_t,\phi_t})_{t \ge 0}$ do not form a nontrivial $C_0$-semigroup on spaces $X$ satisfying 
$H^\infty \subset X \subset \BB$, where $\BB$ is the Bloch space. This includes the case $X=\BMOA$. The proof is based on estimates for derivatives of interpolating Blaschke products.

For composition operators, and with $X=H^\infty$ and $\BB$, this result was shown earlier by Blasco et al \cite{blasco} with an argument 
involving the Dunford--Pettis property. For spaces between $H^\infty$ and $\BB$, the result for composition
operators was given 
by Anderson, Jovovic and Smith \cite{AJS17} using geometric function theory.

\section{Motivation}

\subsection{Universal operators}

Rota \cite{rota59,rota60} introduced the concept of a universal operator.

\begin{defn}
An operator $U \in \LL(H)$ is universal if for all nonzero $T \in \LL(H)$ there is a closed subspace $\MM \ne \{0\}$ of $H$,
an isomorphism $J: \MM \to H$ and a $\lambda \in \CC \setminus \{0\}$
such that $U \MM \subset \MM$ and $U_{| \MM} = J^{-1} (\lambda T) J$.
\end{defn}

That is, a universal operator is a ``model'' for all $T \in \LL(H)$. 
Universal operators are of interest in the study of the {\em invariant subspace problem}, whether
every operator on a separable infinite-dimensional Hilbert space has a nontrivial closed invariant subspace.
This has a positive solution if and only every minimal invariant subspace of a given universal operator is finite-dimensional.
We refer to \cite{CPbook} for more details and   examples.

Caradus \cite{caradus} gives a convenient sufficient condition for
a Hilbert space operator to be universal.

\begin{thm}
Let $U \in \LL(H)$ be such that:\\
(i) $\dim \ker U = \infty$; and\\
(ii) $U$ is surjective.\\
Then $U$ is universal.
\end{thm}

Consider now the hyperbolic automorphisms
\[
\phi_r(z)= \frac{z+r}{1+rz}
\]
with $r \in (-1,1) $. It is helpful to
consider them using the parametrization
\[
\psi_t(z):=\phi_r(z) \qquad \hbox{with} \quad r=\frac{1-e^{-t}}{1+e^{-t}}
\]
for $t \in \RR$.  

It was shown by Nordgren, Rosenthal and Wintrobe \cite{NRW87} that for $r \ne 0$ 
the operator $C_{\phi_r}-\Id$ is
universal on $H^2$. Of course $C_{\phi_r}$ and $C_{\phi_r}-\Id$ have the same invariant subspaces.

A simpler proof using the fact that $C_{\phi_r}$ can be embedded in a $C_0$-group was
given by 
Cowen and Gallardo-Guti\'errez \cite{CGG17}.

Some work on the invariant subspaces of such operators is due to Matache \cite{mat93}, Mortini \cite{mor95}
and Gallardo-Guti\'{e}rrez  and Gorkin \cite{GGG11}.

\subsection{Isometries}

Another application of weighted composition operators goes back to Banach \cite{banach}, who showed that
every surjective isometry $F$ of the   space $C(K)$ of continuous complex functions on a compact metric space $K$
has the form
\[
F (f)  = w(f \circ \phi),
\]
where $w \in C(K)$ satisfies $|w|=1$  and $\phi$ is a homeomorphism of $K$.

The Hardy space $H^2(\DD)$ is a Hilbert space, and thus has many (linear) isometries; however
for other $p$ with $1<p<\infty$ there are relatively few, and they are expressible as
weighted composition operators. In \cite{DLRW} deLeeuw, Rudin, and Wermer
gave a description of the isometric  surjections of $H^1(\DD)$, which arise from
conformal mappings of $\DD$ onto $\DD$.

Moreover, Forelli \cite{forelli} gave the following theorem, which does not assume
surjectivity.

\begin{thm}\label{thm:forellithm}
Suppose that $p \ne 2$ and that $T: H^p(\DD) \to H^p(\DD)$ is a linear isometry.
Then there are a non-constant inner function $\phi$ and a function
$F \in H^p(\DD)$ such that 
\beq\label{eq:forelli15}
Tf= F(f \circ \phi) 
\eeq
for $ f \in H^p(\DD)$.
\end{thm}

\subsection{Change of domain}

It is well known that composition operators on the Hardy space $H^2(\CC_+)$ of the right half-plane
are unitarily equivalent to weighted composition operators on $H^2(\DD)$. For example
the following explicit formula is given in \cite{CP03}.

\begin{prop}
Let $M$ denote the self-inverse bijection from $\DD$ onto $\CC_+$ given by $M(z)=\dfrac{1-z}{1+z}$, and let
$\Psi: \CC_+ \to \CC_+$ be holomorphic.
Then the composition operator $C_\Psi$ on $H^2(\CC_+)$ is unitarily equivalent to the
operator $L_\Phi: H^2(\DD) \to H^2(\DD)$ defined by
\[
L_\Phi f(z)= \frac{1+\Phi(z)}{1+z} f(\Phi(z)),
\]
where $\Phi = M \circ \Psi \circ M$.
\end{prop}

So for example the $C_0$-group $(T_t)_{t \in \RR}$ on $H^2(\CC_+)$ given by $T_t g(z)=g(e^t z)$ $(z \in \CC_+)$
for $g \in H^2(\CC_+)$
is unitarily equivalent to the weighted composition group $(S_t)_{t \in \RR}$ given by
\[
S_t f(z) = \frac{2}{1+z+e^t(1-z)} f \left( \frac{1+z-e^t(1-z)}{1+z+e^t(1-z)}\right).
\]

Formulae for general domains are given, for example, in \cite[Prop. 2.1]{kp04}. If $W$ is a weighted composition operator
between two Hardy--Smirnoff spaces $E^2(\Omega_1)$ and $E^2(\Omega_2)$, with $\Omega_1$ and 
$\Omega_2$ conformally equivalent to the disc $\DD$, then $W$ is unitarily equivalent to
a weighted composition operator on $H^2(\DD)$. Similar formulae are given for Bergman spaces.

\section{Asymptotic behaviour of $T^n$ or $T_t$}
\label{sec:4ab}

\subsection{The discrete unweighted case}\label{sec:3.1}
\label{sec:duc}

For a fixed composition operator $C_\phi$ there are several possible modes of convergence for the sequence $(C^n_\phi)_{n \ge 1}$, some of which we now list in progressively weaker order.
\begin{itemize}
\item Norm convergence. There exists an operator $P \in \LL(X)$ such that $\|C_\phi^n - P \| \to 0$.
\item Strong convergence. There exists an operator $P \in \LL(X)$ such that $\|C_\phi^nx - P x\| \to 0$
for all $x \in X$.
\item Weak convergence. There exists an operator $P \in \LL(X)$ such that $C_\phi^nx  \to P x$ weakly for all $x \in X$.
\end{itemize}
In each case $P$ is the projection onto $\fix(C_\phi):= \{x \in X: C_\phi x=x\}$ along the decomposition
\[
X= \fix(C_\phi) \oplus \overline{\Im(\Id-C_\phi)}.
\]

The following theorem from \cite{arendt-nagel} helps with the analysis.

\begin{thm}\label{thm:nagel}
Let $T \in \LL(X)$ with $\sup_n \|T^n\|< \infty$. Then the following are equivalent.\\
(i) $P:=\lim T^n$ exists and $P$ is a finite-rank operator.\\
(ii) (a) The essential spectral radius $r_e(T)$ satisfies $r_e(T)<1$;\\
\  (b) $\sigma_p(T) \cap \TT \subset \{1\}$; and\\
\  (c) if $1 \in \sigma(T)$ then $1$ is a pole of the resolvent $(z\Id-T)^{-1}$ of order 1.
\end{thm}
In this case $P$ is the residue at 1.\\

We sketch the proof.

\beginpf
(i) $\implies$ (ii):\\
Let $X_1=PX$, $X_2=(\Id-P)X$ and $T_i=T_{X_i}$ for $i=1,2$.

Then $\|T_2^n\|_{\LL(X_2)} \to 0$ as $n \to \infty$, and hence $r(T_2)<1$.

Since $\sigma(Y)=\sigma(T_1) \cup \sigma(T_2) = \{1\} \cup \sigma(T_2)$
and \[
(\lambda \Id-T)^{-1}=
\begin{cases}
\frac{1}{\lambda-1} & \hbox{on } X_1, \\
(\lambda \Id-T_2)^{-1} & \hbox{on }X_2,
\end{cases}
\]
we see that (ii) follows.\\

\noindent (ii) $\implies$ (i):\\
Let $P$ be the residue at $1$, and let $X_1=PX$, $X_2=(\Id-P)X$ and $T_i=T_{X_i}$ for $i=1,2$.

Then $\sigma(T_1)=\{1\}$ and $\sigma(T_2)=\sigma(T) \setminus \{1\}$ by (a) and (b).

Thus $r(T_2)<1$ and so $\|T_2^n\|_{\LL(X_2)} \to 0$.

It follows from (c) that $T_1$ is diagonalisable and thus $T_1=\Id$.
\endpf

Recall that $X \hookrightarrow \Hol(\DD)$ means that $X$ embeds continuously in $\Hol(\DD)$, which 
means that, for all $\lambda\in\DD$,  $\delta_\lambda: f\mapsto f(\lambda)$ from $X$ to $\CC$ is bounded. \\

Arendt and Batty \cite{AB88} have given criteria for strong convergence. More recently,
from \cite{ACKS18} we mention the following theorem.

\begin{thm}
Let $\phi: \DD \to \DD$ be holomorphic,  $C_\phi \in \LL(X)$, where $X$ is a Banach space such that  
$X \hookrightarrow \Hol(\DD)$. Then $(C_\phi^n)_n$ converges uniformly if and only if $r_e(C_\phi)<1$.
\end{thm}

So it remains to study the essential spectral radius of $C_\phi$. To this end we write $\phi^{(n)}=\underbrace{\phi \circ \ldots \phi}_{  \hbox{ $n$ factors}}$.

\subsubsection*{Denjoy--Wolff theory}\label{sec:DW}

A mapping $\phi: \DD \to \DD$ is an {\em elliptic automorphism\/} if it has the form
\[
\phi = \psi_a \circ R_\theta \circ \psi_a,
\]
where $a \in \DD$ and $\psi_a(z)=\dfrac{a-z}{1-\overline a z} = \psi_a^{-1}(z)$,
and 
$R_\theta(z)=e^{i\theta }z$ with $\theta \in \RR$.

Then the classical Denjoy--Wolff theorem states that provided $\phi$ is not an elliptic automorphism, the
sequence $(\phi^{(n)})_n$ converges uniformly to some $\alpha \in \overline\DD$ on each compact subset of $\DD$.

Such an $\alpha$ is called the {\em Denjoy--Wolff point\/} of $\phi$ and is sometimes denoted by $\DW(\phi)$.\\

{\bf Case 1.} $|\alpha|=1$.

\begin{thm} Suppose that $\CC[z] \subset X$ and suppose that each $L \in X'$ with $L(e_n)=L(e_1)^n$ has the form
$L(f)=\delta_{z_0}(f):=f(z_0)$ for some $z_0 \in \DD$. Then $\sup_n \|C_\phi^n\|=\infty$.
Therefore even weak convergence of $(C_\phi^n)$ cannot occur.
\end{thm}

\beginpf
If $\sup_n \|C_\phi^n\|=M < \infty$ then
\[
|f(\alpha)|=\lim_{n \to \infty}|f(\phi^{(n)}(0))| \le \|\delta_0\| M \|f\|,
\]
so $\delta_\alpha=\delta_{z_0}$ for some $z_0 \in \DD$, which is absurd.
\endpf

{\bf Case 2.} $|\alpha|<1$.\\

The natural candidate for $P$ is given by $Pf=f(\alpha)=\delta_\alpha(f){\bf 1}$, a rank-one operator.

\begin{thm}
For composition operators on $H^p(\DD)$, with $1 \le p < \infty$, the
sequence $(C_\phi^n)$ converges uniformly and strongly if and only if $|\alpha|<1$ and $\phi$ is not inner.
It converges weakly if and only if $|\alpha|<1$.
\end{thm}

We proved that $r_e(C_\phi)<1$ on $H^2$ if and only if
$\phi$ is not inner and $|\alpha|<1$.
When $\phi$ is inner and $|\alpha|<1$ then $C_\phi$
is similar to an isometry (this is a necessary and sufficient condition \cite{bayart}).
We deduce the result for $H^p(\DD)$ from a theorem of Shapiro \cite{shap87}, namely
$r_{e, H^2}(C_\phi))^2=r_{e, H^p}(C_\phi))^p$. For $p=1$ there is the inequality ``$>$''.
Hence $r_{e, H^2}(C_\phi))<1$ implies that $r_{e, H^p}(C_\phi))<1$.

We may summarise the results obtained on various Banach spaces in  a table \cite{ACKS18, ACKS20}.

\begin{center}
\begin{tabular}{ | m{5.5cm} | m{1.5cm}|m{1.5cm} | m{1.5cm} |} 
\hline
Space & Uniform  & Strong  & Weak \\
  \hline
  $H^p(\DD)$ & $\phi$~not inner, $|\alpha|<1$  & $\phi$~not inner, $|\alpha|<1$ &  $|\alpha|<1$ \\ 
  \hline
 $A^p_\beta(\DD)$, $\beta>-1$, $1 \le p < \infty$ & $|\alpha|<1$ & $|\alpha|<1$ & $|\alpha|<1$ \\ 
  \hline
 $\BB_0, \BB^\gamma$  & $|\alpha|<1$ & $|\alpha|<1$ & $|\alpha|<1$ \\ 
  \hline
 $H_{\nu_p}^\infty$, $0<p<\infty$ &  $|\alpha|<1$ & $|\alpha|<1$ & $|\alpha|<1$ \\
  \hline
\end{tabular}
\end{center}
We write  $dA$ for the normalized Lebesgue area measure on  $\DD$, i.e.  $dA(re^{i\theta})=\frac{1}{\pi}rdrd\theta$. The \emph{standard weighted Bergman space}, $A_\beta^p(\DD)$,   $\beta\geq -1$, $p\geq 1$ is the space of all holomorphic functions $f:\DD\to\CC$ such that 
\[
\int_{\DD}|f(z)|^p(1-|z|^2)^\beta dA(z) < \infty.
\]
Every $A_\beta^p$ is a Banach space when $1\leq p<\infty$ with norm the $p^{th}$ root of above integral, denoted by $\|f\|_{A_\beta^p}$. \\
The \emph{unweighted Bergman space}, $A^p$ is obtained when $\beta=0$.\\
The \emph{standard Hardy space} $H^p(\DD)$ are obtained when $\beta=-1$.

Here $\BB_0$ and $\BB_\gamma$ are variations on the Bloch space, namely,
\[
\BB_0 = \{ f \in \Hol(\DD): \lim_{|z| \to 1} (1-|z|^2) |f'(z)| = 0 \}
\]
and
\[
\BB^\gamma = \{ f \in \Hol(\DD): \sup_{z \in \DD} (1-|z|^2)^\gamma |f'(z)| <\infty \}.
\]
For $p>0$, the \emph{standard weighted Bergman space of infinite order}, $H_{\nu_p}^\infty(\DD)$ (or $H_{\nu_p}^\infty$), is the  Banach space of all holomorphic functions $f:\DD\to\CC$ such that 
\[
\|f\|_{H_{\nu_p}^\infty}:=\sup_{z\in\DD}\nu_p(z)|f(z)| < \infty,
\]
with the norm as defined above,  where $\nu_p(z)= (1-|z|^2)^p$. 

\subsection{The continuous  unweighted case}\label{sec:3.2}

We now address the question whether we can deduce the asymptotic behaviour of a $C_0$
semigroup $(C_{\phi_t})_{t \ge 0}$
from properties of its generator.

As we shall explain in more detail below, for  $X \hookrightarrow \Hol(\DD)$ with $(C_{\phi_t})_{t \ge 0}$ a $C_0$-semigroup induced by
an analytic semiflow $( \phi_t)_{t \ge 0}$ with 
\[
G(z)=\lim_{t \to 0} \dfrac{\phi_t(z)-z}{t},
\] 
the generator $A$ is given by $Af=Gf'$ ($f \in D(A)$) with dense domain $D(A)=\{f \in X: Gf' \in X \}$.

Various well-known properties of analytic semiflows are the following, which can be found in
the recent book \cite{BCD20}.

\begin{enumerate}
\item For all $t \ge 0$, $\phi_t$ is injective.
\item If there is a $t_0>0$ such that $\phi_{t_0}$ is an (elliptic) automorphism,
then $\phi_t$ is an elliptic automorphism for all $t > 0$.
\item For all semiflows that are not elliptic automorphisms, there is a unique $\alpha \in \overline \DD$
such that $\lim_{t \to 0} \|\phi_t-\alpha\|_{\infty,K}=0$ for all compact subsets $K \subset \DD$. Such an 
$\alpha$ is called the Denjoy--Wolff point of $( \phi_t)_{t \ge 0}$ (see also Section \ref{sec:DW}).
\end{enumerate}

The following classical theorem of Berkson and Porta \cite{BP78}
describes a semigroup in terms of its generator.

\begin{thm}\label{thm:BP78}
let $(\phi_t)_{t \ge 0}$ be an analytic semiflow on $\DD$. Then the generator
\[
G(z):= \lim_{t \to 0} \dfrac{\phi_t(z)-z}{t}
\]
exists for all $z \in \DD$. Also $G \in \Hol(\DD)$ and
\beq\label{eq:Ggen}
G(z)=(\alpha-z)(1-\overline\alpha z)F(z)
\eeq
where $F \in \Hol(\DD)$ and $\re F \ge 0$ (this implies that $F \in \bigcap_{0<p<1} H^p(\DD)$ and hence has radial
limits almost everywhere on $\TT$). Reciprocally, any $G$ of the form \eqref{eq:Ggen} is the generator of a semiflow.
\end{thm}

Another characterization was given in \cite[Thm. 3.9]{ACP15}. If one knows a priori that    $G\in H^1(\DD)$, the necessary and sufficient condition is that
  $\re (\overline z G^*(z)) \le 0$ a.e.\ on $\TT$, where $G^*$ denotes the radial limit of $G$.  
  The proof relies on the following  observation which can be seen as  a maximum principle for  harmonic functions.  If $h:\DD\to\CC$ is in $H^1(\DD)$, then $h(z)=\frac{1}{2\pi}\int_0^{2\pi} \re \left( \frac{e^{it} +z}{e^{it}-z}\right)  f^*(e^{it})dt$, where $h^*$ stand for the radial limit of $h$,  and thefore if $\re (h^*(e^{it}))\leq 0$ a.e. on $\TT$, then $\re(h(z))\leq 0$ for all $z\in \DD$.       \\

We have already mentioned that $F$ and thus $G$ are always in $H^p(\DD)$ for $0<p<1$, which implies the existence of $G^*$ in $L^p(\TT)$. Nevertheless  it is not possible to improve \cite[Thm. 3.9]{ACP15}  since $G(z):=-z\left( \frac{z-1}{z+1}\right)$ is in $H^p(\DD)$ for all $p<1$, satisfies $\re (\overline z G^*(z)) \le 0$ a.e.\ on $\TT$  and is not the generator of  a semiflow using \eqref{eq:Ggen}.   \\

Combining \cite[Thm. 10.2.6]{BCD20} which is called ``Abate's formula" and   \cite[Thm. 10.2.10]{BCD20}, we obtain the following characterization of generators. 
\begin{thm}\label{thm:BCD20}
	Let $G:\DD\to\CC$ be  holomorphic function.  The following assertions are equivalent:
	\begin{itemize}
		\item[(i)] $G$  is the generator of an  analytic semiflow on $\DD$;
		\item[(ii)] $\re (2\overline{z}G(z) +(1-|z|^2)G'(z))\leq 0$ for all $z\in\DD$;
		\item[(iii)] $\limsup_{z\in\DD,z\to\xi} \re(\overline{z}G(z))\leq 0$ for all $\xi\in\TT$.   
	\end{itemize}  
\end{thm}

Now, from results on $H^p(\DD)$ and similar spaces  in Section \ref{sec:DW} we obtain the following result.

\begin{thm}
Let $X=H^p(\DD)$ or other classical Banach spaces (e.g. Berg\-man, Dirichlet). Let $(C_{\phi_t})_{t \ge 0}$
be a $C_0$-semigroup of composition operators on $X$ with generator $A:f \mapsto Gf'$. Then the semigroup
converges uniformly if and only if it converges strongly, and this is if and only if
$G$ has a zero in $\DD$, and there exists $\xi\in\TT$  
such that $\limsup_{z\in\DD,z\to\xi}\re (\overline z G(z)) < 0$.  
\end{thm}

Indeed, $G$ has a zero in $\DD$ if and only if either $\phi_t$ is an elliptic automorphism for each $t$ or else
the Denjoy--Wolff point $\alpha$ lies in $\DD$. The condition $\limsup_{z\in\DD,z\to\xi}\re (\overline z G(z)) < 0$  
implies that $(\phi_t)_{t \ge 0}$ is not a semiflow of automorphisms. 
For if  $(\phi_t)_{t \ge 0}$ were a semigroup of automorphisms, then both $G$ and $-G$ would
generate semigroups and so we would have $\lim_{z\in\DD,z \to \xi}\re (\overline z G(z)) = 0$ everywhere on $\TT$.
Note that since $\phi_t$ is injective, it is not inner precisely when it is not an automorphism.

\subsection{Weighted composition operators}

Iterates of weighted composition operators on holomorphic function spaces $X$ that embed continuously
in $\Hol(\DD)$ are studied in \cite{CP21}.
Some of the most relevant theorems here are the following, where again $W_{w,\phi}$ is the
weighted composition operator $f \mapsto w(f \circ \phi)$. The methods of proof are similar to those of \cite{ACKS18}.
We assume first that $\phi$ is not an elliptic automorphism and that $W_{w,\phi}$ is at least power bounded.

\begin{thm}\label{thm:2.2}
With $X$ and $W_{w,\phi}$ as above, suppose also that $\alpha:=\DW(\phi) \in \DD$. Then
 the sequence $(W_{w,\phi}^n)$ converges weakly as $n \to \infty$
 if
and only if (i) $|w(\alpha)| <1$, or (ii) $|w(\alpha)| =1$ and $\sup_n \|W_{w,\phi}^n\|< \infty$.
\end{thm}

\begin{thm}
Under the hypotheses of Theorem \ref{thm:2.2}, suppose that
either $|w(\alpha)| <1$, or (ii) $|w(\alpha)| =1$ and $\sup_n \|W_{w,\phi}^n\|< \infty$. Then
 $(W_{w,\phi}^n)$ converges  uniformly if and only if $r_e(W_{w,\phi}) < 1$.
 \end{thm}
 
 For elliptic automorphisms the story is slightly different.
 
\begin{thm}
Let $\phi$ be an elliptic automorphism of infinite order with
fixed point $\alpha\in \DD$ and $w \in A(\DD)$ bounded away from 0 on $\DD$. Suppose
that $\sup_n \|W_{w,\phi}^n\|< \infty$. Then $(W_{w,\phi}^n)$ converges  uniformly if and only if
$|w(\alpha)| < 1$.
\end{thm}

\subsection{Isometry and similarity to isometry}

As a byproduct of the results in Section \ref{sec:duc} we have
characterizations of the properties isometry and similarity to an isometry
for composition operators \cite{ACKS18,ACKS20}.
\begin{center}	
	{{
			\begin{tabular}{ | m{5.5cm} | m{3.5cm}| m{3.5cm} |  } 
				\hline 
				\textbf{Spaces}	& $C_\varphi$ \textbf{isometric}&$C_\varphi$ \textbf{similar to an isometry} \\ 
				\hline 
				$H^p(\DD),1\leq p<\infty$	&  $\varphi$ inner and $\varphi(0)=0$  &$\varphi$ inner and $\exists b\in\DD$ with $\varphi(b)=b$   \\ 
				\hline 
				$A^p_\beta(\DD),\beta>-1,1\leq p<\infty$		& $\varphi$ rotation &$\varphi$ elliptic automorphism\\ 
				\hline 
				$\BB$	& $\varphi(0)=0$ and $\tau_{\varphi}^\infty=1$ & $\exists b\in\DD$ with $\varphi(b)=b$  and $\tau_{\varphi}^\infty=1$  \\ 
				\hline 
				$\BB_0$ or $\BB^\alpha$, $\alpha\neq 1$	&  $\varphi$ rotation &  $\varphi$ elliptic automorphism\\ 
				\hline 
				$H^\infty_{\nu_p},0<p<\infty$	& $\varphi$ rotation & $\varphi$ elliptic automorphism \\ 
				\hline 
			\end{tabular}   	
	}  } \\
\end{center} 


Here 
\[
\tau_\phi^\infty : = \sup_{z \in \DD} \frac{1-|z|^2}{1-|\phi(z)|^2} |\phi'(z)|.
\]

\subsection{Generators}


Suppose now that $(T_t)_{t \ge 0}$ is a $C_0$-semigroup of weighted composition operators
$T_t = w_t C_{\phi_t}$ acting on a suitable space of functions on a domain $\Omega$. We have already discussed
results 
describing generators in the unweighted case
in Theorems \ref{thm:BP78} and \ref{thm:BCD20},
namely $Af=Gf'$, where
\beq \label{eq:Gz}
G(z):= \lim_{t \to 0} \dfrac{\phi_t(z)-z}{t}.
\eeq
The converse, as shown in work by Arendt and Chalendar \cite{AC19} and
Gallardo and Yakubovich \cite{GY19}, associates with a holomorphic function $G$ a semiflow
satisfying $u'(t)=G(u(t))$, with $u(0)=z$. 

\begin{thm}\label{thm:acgy}
Let $X$ be a function space with continuous point evaluations on a domain $\Omega$. If  either \\
(i)  If $(z_n) \subset \Omega$ such that $z_n \to z \in  \Omega \cup \{\infty\}$ and $ \lim_{n \to \infty}
f(z_n)$ exists in $\CC$ for all 
$f \in X$ then $z \in \Omega$, or\\
(ii) $\Omega=\DD$ and $\Hol(\overline \DD) \subset X \subset \Hol(\DD)$ with continuous embeddings,\\
then the semigroup associated with $G$ is a (quasicontractive) semigroup of composition operators.
\end{thm}

Note that $\Hol(\overline \DD)$ embeds continuously in  $X$ means that for each $\varepsilon>0$,  there exists a positive constant $c(\varepsilon)$ such that 
\[\|f\|_X\leq c(\epsilon) \sup_{|z|\leq 1+\epsilon}|f(z)| ,\,\, f\in A(D(0,1+\varepsilon)),\]
where   $A(D(0,1+\varepsilon))$ is the algebra of all holomorphic functions in $D(0,1+\varepsilon)$, the open disc centered at $0$ and of radius $1+\varepsilon$, and continuous on  the closure of   $D(0,1+\varepsilon)$. \\

For weighted composition operators on $\DD$, a semigroup of the form $T_t = w_t C_{\phi_t}$ has infinitesimal
generator $Af=Gf'+gf$, where now 
\[
g(z)=\frac{\partial w_t(z)}{\partial t}_{t=0}
\]
and $G$ is once again given by \eqref{eq:Gz}. Assuming Condition (ii) in Theorem \ref{thm:acgy}, it is shown in \cite[Thm. 3.1]{GSY21} that every semigroup with generator of the
form $Af=Gf'+gf$ is a semigroup of weighted composition operators. See also \cite{bernard} for a generalization of \cite{AC19}.



\section{Compactness and analyticity}
This section is mainly extracted from \cite{ACP16}.  We restrict our study to the  semigroups of composition operators on the Hardy space $H^2(\DD)$, even though some of the results can be extended to non-Hilbertian Hardy spaces or weighted Hardy spaces such as the Dirichlet space.     
\subsection{Immediate and eventual compactness}
We recall that a semigroup $(T(t))_{t \ge 0}$ is said to be {\em immediately compact}    
if
the operators $T(t)$ are compact   for all $t>0$. 
A semigroup $(T(t))_{t \ge 0}$ is said to be {\em eventually compact} if there exists $t_0>0$ such that $T(t)$ is compact for all $t\geq t_0$.  
Similar definitions hold for immediately/eventually   trace-class.

The following theorem \cite[Chap. 2, Thm 3.3]{pazy}  links
immediate compactness with continuity in norm.
\begin{thm}\label{thm:cpresolv}
	Let $(T(t))_{t \ge 0}$ be a $C_0$-semigroup and let $A$ be its infinitesimal generator.
	Then $(T(t))_{t \ge 0}$ is immediately compact if and only if\\ (i)~the resolvent $R(\lambda,A)$ is
	compact for all (or for one) $\lambda \in \CC \setminus \sigma(A)$, and\\ 
	(ii)~$\lim_{s \to t} \|T(s)-T(t)\| =0$ for all $t>0$.
\end{thm}

We begin with an elementary observation.

\begin{prop}\label{prop:notimcomp}
	Suppose that for some $t_0>0$ one has $|\phi_{t_0}(\xi)|=1$ on a set of positive measure; then
	$C_{\phi_{t_0}}$ is not compact on $H^2(\DD)$  and so
	the semigroup $(C_{\phi_t})_{t \ge 0}$ is  not immediately compact.
\end{prop}

\beginpf
For the Hardy space, this follows since the weakly null sequence
$(e_n)_{n \ge 0}$ with $e_n(z)=z^n$ is mapped into $(\phi_{t_0}^n)$, which does not converge to
$0$ in norm.
\endpf

We shall now give a sufficient condition for immediate compactness of a semigroup of composition
operators, in terms of the associated function $G$. First, we recall a classical necessary and
sufficient condition for compactness of a composition operator $C_\phi$  in the case when
$\phi$ is univalent \cite[pp. 132, 139]{CM95}.

\begin{thm}
	For $\phi: \DD \to \DD$ analytic and univalent, the composition operator $C_\phi$ is compact
	on $H^2(\DD)$
	if and only if
	\[
	\lim_{z \to \xi} \frac{1-|z|^2} {1-|\phi(z)|^2} = 0
	\]
	for all $\xi \in \TT$.
\end{thm}
The following proposition collects together standard results on  
trace-class composition operators \cite[p.~149]{CM95}.

\begin{prop}\label{prop:hs}
	For $\phi:\DD\to\DD$ analytic  with $\|\phi\|_\infty<1$, the composition operator $C_\phi$ is 
	trace-class on $H^2(\DD)$. 
\end{prop}

Here is an example showing that immediate and eventual compactness are different even for semigroups of composition operators  on $H^2(\DD)$.

\begin{ex}
	Let $h$ be the Riemann map from $\DD$ onto the starlike region 
	\[\Omega:=\DD\cup \{z\in\CC:0<\re (z) <2\mbox{ and }0<\Im(z)<1\},\]
	with $h(0)=0$. Since $\partial\Omega$ is a Jordan curve, 
	the Carath\'eodory theorem \cite[Thm 2.6, p. 24]{pomm}
	implies that $h$ extends continuously to $\TT$.
	 
	Let $\phi_t(z)=h^{-1}(e^{-t}h(z))$.  Note that for $0<t<\log 2$, $\phi_t(\TT)$ intersects $\TT$ on a  set of positive measure, and thus, $C_{\phi_t}$ is not compact
	by Proposition \ref{prop:notimcomp}.
	Moreover, for $t>\log 2$, $\|\phi_t\|_\infty<1$, and therefore $C_{\phi_t}$ is compact (actually trace-class). Figure 1 represents the image of $\varphi_t$  for different values of $t$.  
	
\begin{center}
		\includegraphics[width=14.5cm, trim = 15mm 200mm 0.5cm 0, clip  ]{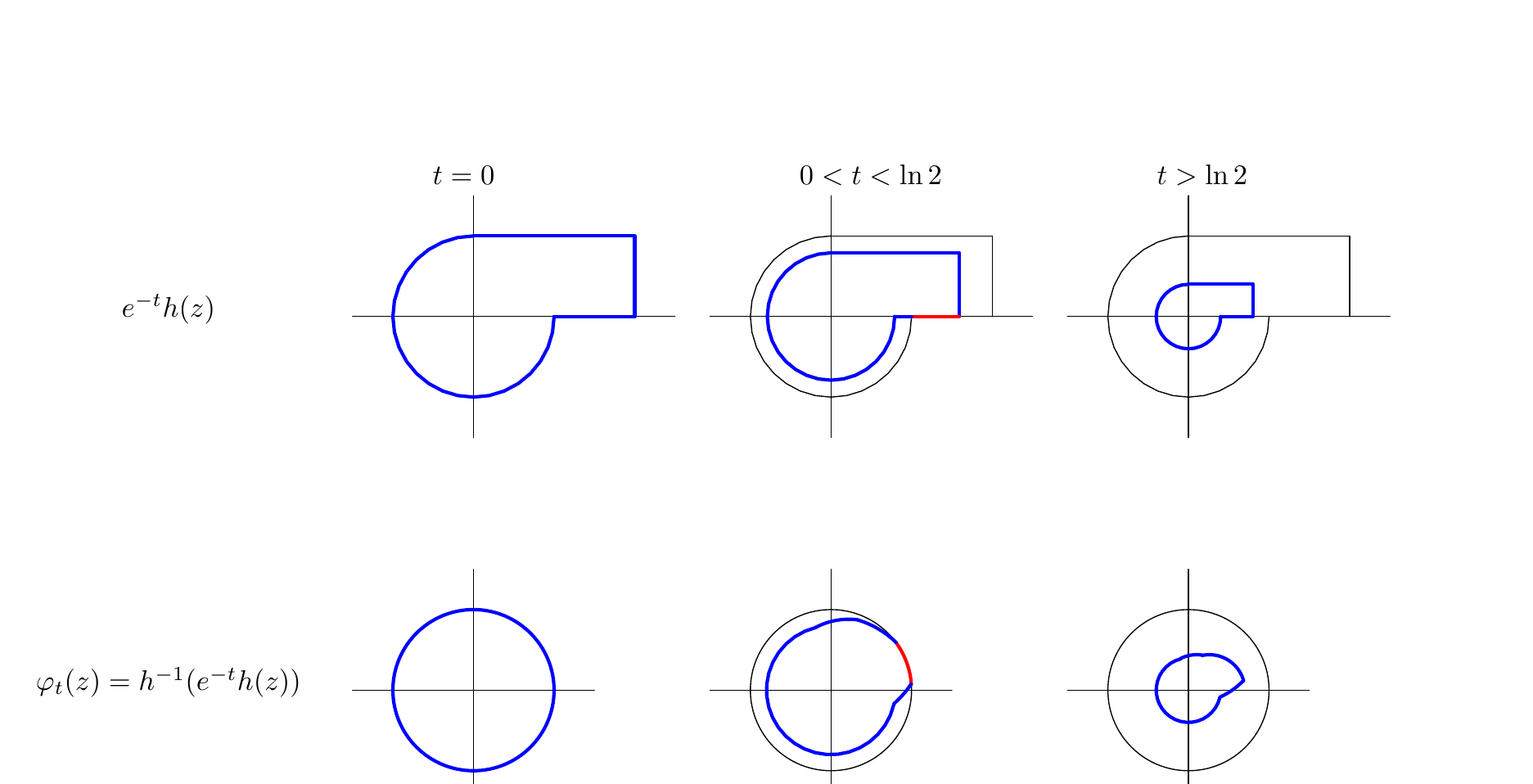}
		Figure 1.
	\end{center}
\end{ex}

\subsection{Compact analytic semigroups}

A $C_0$-semigroup $T$ will be called analytic (or holomorphic) if there exists a sector $\Sigma_\theta=\{ re^{i\alpha}, r\in\RR_+, |\alpha|<\theta\}$ with $\theta \in (0,\frac\pi 2]$ and an analytic mapping $\widetilde T:\Sigma_\theta \to \mathcal L(X)$ such that $\widetilde T$ is an extension of $T$ and
\[\sup_{\xi\in \Sigma_\theta \cap \DD} \| \widetilde T (\xi)\| <\infty.\]
In both cases, the generator of $T$ (or $\widetilde T$) will be the linear operator $A$ defined by
\[
D(A)=\left\{x\in X, \lim_{\RR\ni t\to 0} \frac{T(t)x-x}{t} \text{ exists}\right\}\]
and, for all $x \in D(A)$,
\[  Ax=\lim_{\RR\ni t\to 0} \frac{T(t)x-x}{t}.\]

 In the particular case of analytic semigroups, the compactness is equivalent to the compactness of the resolvent,
 by Theorem \ref{thm:cpresolv}, since the analyticity implies the uniform continuity \cite[p. 109]{EN06}. 
 
 \begin{rem}\label{rem:pisier}
 		For an analytic semigroup $(T(t))_{t\geq 0}$, being eventually compact is equivalent to being immediately compact.  Indeed, consider $Q$, the quotient map from the  bounded  linear operators on a Hilbert space ${\mathcal L}(H)$ onto the Calkin algebra (the quotient of ${\mathcal L}(H)$ by the compact operators). Then $(QT(t))_{t\geq 0}$ is an analytic semigroup which vanishes for $t>0$ large enough, and therefore vanishes identically \
 		(this observation is attributed to W.~Arendt).   
 \end{rem}
 
 We may include that remark in the following result.
 
 \begin{thm}\label{th:compact}
	Let $G: \DD \to \CC$ be a holomorphic function such that the operator  $A$ 
	defined by $Af(z)=G(z)f'(z)$ with dense domain $D(A)\subset H^2(\DD)$  
	generates an analytic semigroup $(T(t))_{t\geq 0}$ of composition operators.  Then the following assertions are equivalent:
	\begin{enumerate}
		\item $(T(t))_{t\geq 0}$ is immediately compact;
		\item $(T(t))_{t\geq 0}$ is eventually compact;
		\item 
		$\forall \xi\in\TT$, $\lim_{z\in\DD,z\to\xi }\left| \frac{G(z)}{z-\xi}\right|=\infty$. 
	\end{enumerate} 
\end{thm}

\begin{thm}\label{th:analytic}
	Let $G: \DD \to \CC$ be a holomorphic function such that the operator  $A$ 
	defined by $Af(z)=G(z)f'(z)$ has dense domain $D(A)\subset H^2(\DD)$. 
	The operator $A$ generates an analytic semigroup of composition operators on $H^2(\DD)$  if and only if there exists $ \theta \in (0,\frac \pi 2)$ such that for all $\alpha \in \{-\theta,0,\theta\}$
	\[\limsup_{z\in\DD,z\to \xi} \re(e^{i\alpha}\overline z G(z)) \le  0\mbox{ for all }\xi\in\TT.\]
\end{thm}

Using the semiflow model, we have the following result.

\begin{thm}
	Let $(C_{\phi_t})_{t \ge 0}$  be an immediately compact analytic semigroup on $H^2(\DD)$.  
	Then the following conditions are equivalent:\\
	1. There exists a $t_0 >0$ such that $\|\phi_{t_0}\|_\infty < 1$;\\
	2. For all $t>0$ one has $\|\phi_t\|_\infty < 1$.\\
	Therefore, if  there exists a $t_0 >0$ such that $\|\phi_{t_0}\|_\infty < 1$, then
	$(C_{\phi_t})_{t \ge 0}$
	is immediately trace-class.
\end{thm}
It is of interest to consider the relation between immediate compactness and analyticity for
a $C_0$-semigroup of composition operators: this is because compactness of a
semigroup $(T(t))_{t \ge 0}$ is implied by compactness of the resolvent together with norm-continuity
at all points $t>0$, as in Theorem \ref{thm:cpresolv}.

\begin{ex}
	Consider
	\[
	G(z)= \frac{2z}{z-1},
	\]
	Now the image of the unit disc under $\overline z G(z)$ is contained in the left half-plane
	so the operator $A: f \mapsto Gf'$ generates a non-analytic $C_0$-semigroup 
	of composition operators $(C_{\phi_t})_{t \ge 0}$ on $H^2(\DD)$.
	On the other hand, it can be shown that 
	$C_{\phi_t}$ is compact -- even trace-class -- for each $t>0$. For
	we have the equation
	\[
	\phi_t(z)e^{-\phi_t(z)}=e^{-2t}ze^{-z}.
	\]
	Now the function $z \mapsto ze^{-z}$ is injective on $\overline\DD$; this follows from
	the argument principle,  for the image of $\TT$ is easily seen to be a simple Jordan curve.
	It follows that
	$\|\phi_t\|_\infty<1$ for all $t>0$, and so $C_{\phi_t}$ is trace-class.
\end{ex}

\begin{ex}
	The semigroup corresponding to $G(z)=(1-z)^2$ is analytic but not immediately compact.
	For 
	\[
	\phi_t(z)=\frac{(1-t)z+t}{-tz+1+t}
	\]
	 the Denjoy--Wolff point is $1$, so the semigroup cannot be immediately compact.
	
	The analyticity follows on calculating $\overline z G(z)$ for $z=e^{i\theta}$.
	We obtain $-4 \sin^2(\theta/2)$, which gives the result by Theorem~\ref{th:analytic}.
\end{ex}

\section{An outlook on $\CC_+$ and $\CC$}

\subsection{The right  halfplane $\CC_+$}

Unlike in the case of the disc, not all composition operators on $H^2(\CC+)$ are bounded,
and there are no compact composition operators.

The following theorem was given by Elliott and Jury \cite{EJ12} (see also \cite{Mat08}).
Recall that the angular derivative $\phi'(\infty)$ of a self-map of $\CC_+$ is defined by
\[
\phi'(\infty)=\lim_{z \to \infty} \frac{z}{\phi(z)}
\]
\begin{thm}
Let $\phi: \CC_+ \to \CC_+$ be holomorphic. The composition operator $C_\phi$
is bounded on $H^2(\CC_+)$ if and only if $\phi$ has finite angular derivative $0<\lambda<\infty$
at infinity, in which case $\|C_\phi\|=\sqrt\lambda$. We also have for the essential norm
that $\|C_\phi\|_e=\|C_\phi\|$
so that there are no compact composition operators on $H^2(\CC_+)$.
\end{thm}

Berkson and Porta   \cite{BP78} gave the following criterion for an analytic
funcction $G$ to generate a one-parameter semiflow of analytic mappings
from $\CC_+$ into itself, namely, solutions to the initial value problem
\[
\frac{\partial \phi_t(z)}{\partial t}=G(\phi_t(z)), \qquad \phi_0(z)=z,
\]
namely the condition 
\[
x \frac{\partial(\re G)}{\partial x} \le \re G \qquad \hbox{on} \quad \CC_+,
\]
where as usual $x=\re z$. In this case the semigroup with generator $A: f \mapsto Gf'$
consists of composition operators.

Arvanitidis \cite{arv} showed that a necessary and sufficient condition for these
composition operators to be bounded is that the non-tangential limit
$\delta:=\angle\lim_{z \to \infty} G(z)/z$ exists, in which case $\|C_{\phi_t}\|=e^{-\delta t /2}$
and the semigroup is quasicontractive.

In \cite{ACP16} it was shown that for $G$ holomorphic on $\CC_+$ a necessary condition for  the operator $A: f \mapsto Gf'$ to generate a quasicontractive
semigroup is
\[
\inf_{z \in \CC_+} \frac{\re G(z)}{\re z} > -\infty,
\]
and for a contractive semigroup this infimum is non-negative.


\subsection{The complex plane $\CC$}

For $\CC$ we present material on the Fock space from \cite{CP22}. The Fock space is
arguably one of the most important spaces of entire functions, and it is possible
to give a complete answer to several questions involving (weighted)
composition operators.

For $1 \le \nu< \infty$ the Fock space $\FF^\nu$ is defined to be the space of entire functions $f$ on $\CC$
such that the norm
\[
\|f\|_\nu:= \left( \frac{\nu}{2\pi} \int_\CC |f(z)|^\nu e^{-\nu z^2/2} \, dm(z) \right)^{1/\nu}
\]
is finite. The space $\FF^2$ is a Hilbert space with orthonormal basis $(\tilde e_n)_{n=0}^\infty $
defined by 
\[
\tilde e_n(z)=\frac{z^n}{\sqrt{n!}}.
\]

\subsubsection{Composition operators}

Boundedness of composition operators on $\FF^2$ was characterized in \cite{CMS03} as follows.

\begin{thm}
For an entire function $\phi$ the composition operator $C_\phi$ is bounded on
$\FF^2$ if and only if either $\phi(z)=az+b$ with $|a|<1$ and $b \in \CC$ or $\phi(z)=az$ with $|a|=1$.
In the case where $|a|<1$ the operator $C_\phi$ is compact.
\end{thm}
Thus there are relatively few bounded composition operators on $\FF^2$ and most natural questions can be answered easily.
When we come to weighted composition operators there will be much more to say.

The discussion here is based mostly on \cite{CP22}, although some results may also be found in
\cite{seyoum}, which concentrates on spectra and mean ergodicity. 

Considering now iterates $(C_\phi^n)_n$ we have
\[
C_\phi^n f(z) = \begin{cases}
f \left(a^nz + \frac{1-a^n}{1-a}b \right) & \hbox{if } a \ne 1, \\
f(z) & \hbox{if } a=1.
\end{cases}
\]
We no wish to use Theorem \ref{thm:nagel} again, so we need to know whether $C_\phi$ is
power-bounded. The following formula is due to \cite{CMS03} in the case $\nu=2$, and \cite{dai}
in general.
\[
\|C_\phi\|= \exp\left( \frac14 \frac{|b|^2}{1-|a|^2} \right) \qquad \hbox{if} \quad \phi(z)=az+b \quad \hbox{with} \quad |a|<1.
\]

\begin{thm}
The asymptotics of iteration of bounded composition operators on $\FF^\nu$ are as follows.
\begin{enumerate}
\item 
If $\phi(z)=az$ with $|a|=1$ and $a \ne 1$, then $(C^n_\phi)$ consists of unitary operators
and does not converge even weakly.
\item If $\phi(z)=az+b$ with $|a|<1$ then $(C_\phi^n)$ converges in norm to the operator 
$T:f \mapsto f\left( \dfrac{b}{1-a} \right)$.
\end{enumerate}
\end{thm}
The proof is short, and we include it here.
\beginpf
1. In the first case we have that $C_\phi^n f(z)=f(a^n z)$ and for $f(z)=z$ there is clearly no
convergence.

2. In the second case 
\[
\phi^{(n)}(z)=a^nz+ \frac{1-a^n}{1-a}b=: a_n + b_n z, \hbox{ say}.
\]
and 
\[
\frac{|b_n|^2 }{1-|a_n|^2} \to \frac{|b|^2}{|1-a|^2}
\]
so $C_\phi$ is power bounded. Also $r_e(C_\phi)=0$ since $C_\phi$ is compact.

For the point spectrum, suppose that $\lambda$ is an eigenvalue and $f$ an eigenvector. Then
\[
f(az+b)=\lambda f(z)
\]
and so
$f(a_n + b_n z)=\lambda^n f(z) $, which implies that
\beq
\label{eq:lamconv}
\lambda^n f(z) \to f(b/(1-a)).
\eeq
 This means that $f$ is identically zero if $|\lambda|=1$ and $\lambda \ne 1$.

Finally if $1 \in \sigma(C_\phi)$ then $1$ is an eigenvalue since $C_\phi$ is compact. 
By \eqref{eq:lamconv} $f$ is constant, and $\dim\ker (C_\phi-\Id)=1$. Hence $1$ is a pole
of the resolvent of order at most 1. 
Indeed, assume that the pole order of $1$ is larger than $1$. Looking at the Jordan normal form of $T_0:={C_\varphi}_{|X_0}$, where $X_0=P\FF^\nu$ and $P$ the residue, we see that there exists $f\in\FF^\nu$ such that $(\Id -C_\varphi)f=1_{\CC}$. Evaluating at $z=\frac{b}{1-a}$ we obtain a contradiction. Thus the pole order is $1$. It follows that $P$ is the projection onto $\ker (C_\varphi -\Id)=\CC 1_{\CC}$ along 
$\{f\in \FF^\nu:f(b/(1-a))=0\}$. Thus $Pf=f(b/(1-a))1_{\CC}$.     
This completes the proof that $C^n_\phi \to P$ is norm thanks to Theorem~\ref{thm:nagel}. 

\endpf

Similarly, we can characterise $C_0$-semigroups of bounded composition operators.

\begin{thm}
A $C_0$-semigroup $(C_{\phi_t})_{t \ge 0}$ of bounded composition of $\FF^\nu$ satisfies
one of the following conditions for $t>0$.
\begin{enumerate}
\item $\phi_t(z)=e^{\lambda t}z+ C(e^{\lambda t}-1)$ for some $\lambda \in \CC_-=\{z \in \CC: \re z < 0 \}$ and $C \in \CC$.
\item $\phi_t(z)= e^{\lambda t}z$ for some $\lambda \in i\RR$.
\end{enumerate}
The generator is given by
\[
Af(z)= \lambda (z+C) f'(z),
\]
where $C=0$ in case (2).

Moreover $(C_{\phi_t})_{t \ge 0}$ converges in norm a $t \to \infty$ if and only if $\lambda \in \CC_-$,
in which case the limit $T$ satisfies $Tf=f(-C)$.
\end{thm}

Using a theorem from \cite{AC19} we get:
\begin{thm}
Every $C_0$-semigroup on $\FF^2$ with generator of the form $Af=Gf'$ for some
$G \in \Hol(\CC)$ is a semigroup of composition operators with generator satisfying
\[
\begin{array}{rll}
G(z)&= az+b, & \hbox{with } \re a<0, \quad \hbox{or}\\
G(z)&=az, & \hbox{with } a \in i\RR. 
\end{array}
\]
Moreover, the condition
\beq 
\limsup_{|z| \to \infty} \re \overline z G(z) \le 0
\eeq
is necessary and sufficient for such $G$.
\end{thm}

\subsubsection{Weighted composition operators}

For weighted composition operators $\W$ we have the result of Le \cite{Le}, extended by
Hai and Khoi \cite{HaiKhoi} as follows:

\begin{prop}\label{prop:wcobdd}
The weighted composition operator $\W$ with $w$ not identically zero is bounded on $\FF^\nu$ if and only if both the following
conditions hold:
\begin{itemize}
\item $w \in \FF^\nu$;
\item $M(w,\phi):=\sup_{z \in \CC} |w(z)|^2 e^{(|\phi(z)|^2-|z|^2)} < \infty$.
\end{itemize}
Moreover, in this case $\phi(z)=az+b$ with $|a| \le 1$. If $|a|=1$ then we also
have 
\beq\label{eq:wamod1}
w(z)=w(0)e^{- \overline b az}
\eeq
 for $z \in \CC$.
\end{prop}

See also \cite{CG21} for a recent discussion of these results. Note that from \cite{HaiKhoi}
we have the useful inequality
\[
\sqrt{M(w,\phi)} \le \| \W \| \le \sqrt{M(w,\phi)}/|a|.
\]
Thus for power-boundedness of $\W$, it is possible to give a complete result for $|a|=1$.

\begin{thm}
Suppose that $\phi(z)=az+b$ with $|a|=1$, then $\W$ is power-bounded on
$\FF^\nu$ if and only if $|w(0)| \le e^{-|b|^2/2}$.
\end{thm}

The case $|a|<1$ is apparently unsolved, although this result for non-vanishing $w$ (as
in the case of semigroups) can be found in \cite{CG21}.

\begin{thm} 
For $\phi(z)=az+b$ with $|a|<1$ and 
$w$ nonvanishing the operator $\W$ is bounded on $\GF^\nu$ if and only if
\[
w(z)= e^{p+qz+rz^2}
\]
and either:\\
(i) $|r|< \beta/2$, in which case $\W$ is compact on $\GF^\nu$, or\\
(ii) $|r|=\beta/2$ and, with $t=q+\overline b a$, one has either $t=0$ or else
$t \ne 0$ and $r=-\frac{\beta}2\frac{t^2}{|t|^2}$. In case (ii) $\W$ is not compact on $\GF^\nu$.
\end{thm}

Finally, we can use this to characterise $C_0$-semigroups on the Fock space \cite{CP22}.

\begin{thm}\label{thm:c0sg}
A $C_0$-semigroup $T_t f(z)= w_t(z) f(\phi_t(z))$ ($t \ge 0$) of weighted composition operators on the Fock space  $\GF^\nu$ has one of the
following  two expressions:
\begin{enumerate}  
\item $\phi_t(z)= \exp(\lambda t)z+C (\exp (\lambda t)-1) $ for some $\lambda \in \CC_- $ and $C \in \CC$;
in which case $w_t=e^{p_t+q_t z+r_t z^2}$, where explicit formulae for $p_t$, $q_t$ and $r_t$
can be given.
\item $\phi_t(z)=\exp(\lambda t)z+C (\exp (\lambda t)-1) $ for some $\lambda \in i\RR$ and $C \in \CC$, in which case 
$w_t(z)=w_t(0)\exp(\overline C (\exp (\lambda t)-1)  z)$
and moreover 
$w_t(0)=e^{\mu t}e^{|C|^2 (e^{\lambda t}-1)}$ for some $\mu \in \CC$.
\end{enumerate}
\end{thm}


\begin{thebibliography}{10}

\bibitem{Abate89}
M.~{Abate}.
\newblock {\em {Iteration theory of holomorphic maps on taut manifolds}}.
\newblock Commenda di Rende (Italy): Mediterranean Press, 1989.

\bibitem{AJS17}
A.~Anderson, M.~Jovovic, and W.~Smith.
\newblock Composition semigroups on {BMOA} and {$H^\infty$}.
\newblock {\em J. Math. Anal. Appl.}, 449(1):843--852, 2017.

\bibitem{AB88}
W.~Arendt and C.~J.~K. Batty.
\newblock Tauberian theorems and stability of one-parameter semigroups.
\newblock {\em Trans. Amer. Math. Soc.}, 306(2):837--852, 1988.

\bibitem{AC19}
W.~Arendt and I.~Chalendar.
\newblock Generators of semigroups on {B}anach spaces inducing holomorphic
  semiflows.
\newblock {\em Israel J. Math.}, 229(1):165--179, 2019.

\bibitem{ACKS18}
W.~Arendt, I.~Chalendar, M.~Kumar, and S.~Srivastava.
\newblock Asymptotic behaviour of the powers of composition operators on
  {B}anach spaces of holomorphic functions.
\newblock {\em Indiana Univ. Math. J.}, 67(4):1571--1595, 2018.

\bibitem{ACKS20}
W.~{Arendt}, I.~{Chalendar}, M.~{Kumar}, and S.~{Srivastava}.
\newblock {Powers of composition operators: asymptotic behaviour on Bergman,
  Dirichlet and Bloch spaces}.
\newblock {\em {J. Aust. Math. Soc.}}, 108(3):289--320, 2020.

\bibitem{arendt-nagel}
W.~Arendt, A.~Grabosch, G.~Greiner, U.~Groh, H.~P. Lotz, U.~Moustakas,
  R.~Nagel, F.~Neubrander, and U.~Schlotterbeck.
\newblock {\em One-parameter semigroups of positive operators}, volume 1184 of
  {\em Lecture Notes in Mathematics}.
\newblock Springer-Verlag, Berlin, 1986.

\bibitem{arv}
A.~G. Arvanitidis.
\newblock Semigroups of composition operators on {H}ardy spaces of the
  half-plane.
\newblock {\em Acta Sci. Math. (Szeged)}, 81(1-2):293--308, 2015.

\bibitem{ACP15}
C.~Avicou, I.~Chalendar, and J.~R. Partington.
\newblock A class of quasicontractive semigroups acting on {H}ardy and
  {D}irichlet space.
\newblock {\em J. Evol. Equ.}, 15(3):647--665, 2015.

\bibitem{ACP16}
C.~Avicou, I.~Chalendar, and J.~R. Partington.
\newblock Analyticity and compactness of semigroups of composition operators.
\newblock {\em J. Math. Anal. Appl.}, 437(1):545--560, 2016.

\bibitem{banach}
S.~Banach.
\newblock {\em Theory of linear operations}, volume~38 of {\em North-Holland
  Mathematical Library}.
\newblock North-Holland Publishing Co., Amsterdam, 1987.
\newblock Translated from the French by F. Jellett, With comments by A. Pe\l
  czy\'{n}ski and Cz. Bessaga.

\bibitem{bayart}
F.~Bayart.
\newblock Similarity to an isometry of a composition operator.
\newblock {\em Proc. Amer. Math. Soc.}, 131(6):1789--1791 (electronic), 2003.

\bibitem{BP78}
E.~Berkson and H.~Porta.
\newblock Semigroups of analytic functions and composition operators.
\newblock {\em Michigan Math. J.}, 25(1):101--115, 1978.

\bibitem{bernard}
E.~Bernard.
\newblock Weighted composition semigroups on {B}anach spaces of holomorphic
  functions.
\newblock Preprint, 2022.

\bibitem{blasco}
O.~Blasco, M.D. Contreras, S.~D\'{\i}az-Madrigal, J.~Mart\'{\i}nez,
  M.~Papadimitrakis, and A.G. Siskakis.
\newblock Semigroups of composition operators and integral operators in spaces
  of analytic functions.
\newblock {\em Ann. Acad. Sci. Fenn. Math.}, 38(1):67--89, 2013.

\bibitem{BCD20}
F.~Bracci, M.~D. Contreras, and S.~D\'{\i}az-Madrigal.
\newblock {\em Continuous semigroups of holomorphic self-maps of the unit
  disc}.
\newblock Springer Monographs in Mathematics. Springer, Cham, 2020.

\bibitem{caradus}
S.~R. Caradus.
\newblock Universal operators and invariant subspaces.
\newblock {\em Proc. Amer. Math. Soc.}, 23:526--527, 1969.

\bibitem{CG21}
T.~Carroll and C.~Gilmore.
\newblock Weighted composition operators on {F}ock spaces and their dynamics.
\newblock {\em J. Math. Anal. Appl.}, 502(1):125234, 2021.

\bibitem{CMS03}
B.~J. Carswell, B.~D. MacCluer, and A.~Schuster.
\newblock Composition operators on the {F}ock space.
\newblock {\em Acta Sci. Math. (Szeged)}, 69(3-4):871--887, 2003.

\bibitem{CP03}
I.~Chalendar and J.~R. Partington.
\newblock On the structure of invariant subspaces for isometric composition
  operators on {$H^2(\Bbb D)$} and {$H^2(\Bbb C_+)$}.
\newblock {\em Arch. Math. (Basel)}, 81(2):193--207, 2003.

\bibitem{CPbook}
I.~Chalendar and J.~R. Partington.
\newblock {\em Modern approaches to the invariant-subspace problem}, volume 188
  of {\em Cambridge Tracts in Mathematics}.
\newblock Cambridge University Press, Cambridge, 2011.

\bibitem{CP21}
I.~Chalendar and J.~R. Partington.
\newblock Weighted composition operators: isometries and asymptotic behaviour.
\newblock {\em J. Operator Theory}, 86(1):189--201, 2021.

\bibitem{CP22}
I.~Chalendar and J.R. Partington.
\newblock Weighted composition operators on the {F}ock space: iteration and
  semigroups.
\newblock Preprint, 2021, https://arxiv.org/abs/2106.00427.

\bibitem{CM95}
C.~C. Cowen and B.~D. MacCluer.
\newblock {\em Composition operators on spaces of analytic functions}.
\newblock Studies in Advanced Mathematics. CRC Press, Boca Raton, FL, 1995.

\bibitem{CGG17}
C.C. Cowen and E.A. Gallardo-Guti\'{e}rrez.
\newblock A new proof of a {N}ordgren, {R}osenthal and {W}introbe theorem on
  universal operators.
\newblock In {\em Problems and recent methods in operator theory}, volume 687
  of {\em Contemp. Math.}, pages 97--102. Amer. Math. Soc., Providence, RI,
  2017.

\bibitem{dai}
J.~Dai.
\newblock The norm of composition operators on the {F}ock space.
\newblock {\em Complex Var. Elliptic Equ.}, 64(9):1608--1616, 2019.

\bibitem{DLRW}
K.~de~Leeuw, W.~Rudin, and J.~Wermer.
\newblock The isometries of some function spaces.
\newblock {\em Proc. Amer. Math. Soc.}, 11:694--698, 1960.

\bibitem{EJ12}
S.~Elliott and M.~T. Jury.
\newblock Composition operators on {H}ardy spaces of a half-plane.
\newblock {\em Bull. Lond. Math. Soc.}, 44(3):489--495, 2012.

\bibitem{EN06}
K.-J. Engel and R.~Nagel.
\newblock {\em A short course on operator semigroups}.
\newblock Universitext. Springer, New York, 2006.

\bibitem{forelli}
F.~Forelli.
\newblock The isometries of {$H^{p}$}.
\newblock {\em Canad. J. Math.}, 16:721--728, 1964.

\bibitem{GGG11}
E.A. Gallardo-Guti\'{e}rrez and P.~Gorkin.
\newblock Minimal invariant subspaces for composition operators.
\newblock {\em J. Math. Pures Appl. (9)}, 95(3):245--259, 2011.

\bibitem{GSY21}
E.A. Gallardo-Guti\'{e}rrez, A.G. Siskakis, and D.V. Yakubovich.
\newblock Generators of {$C_0$}-semigroups of weighted composition operators.
\newblock Preprint, 2021, https://arxiv.org/abs/2110.05247.

\bibitem{GY19}
E.A. Gallardo-Guti\'{e}rrez and D.V. Yakubovich.
\newblock On generators of {$C_0$}-semigroups of composition operators.
\newblock {\em Israel J. Math.}, 229(1):487--500, 2019.

\bibitem{HaiKhoi}
P.~V. Hai and L.~H. Khoi.
\newblock Boundedness and compactness of weighted composition operators on
  {F}ock spaces {$\FF^p(\CC)$}.
\newblock {\em Acta Math. Vietnam.}, 41(3):531--537, 2016.

\bibitem{kp04}
R.~Kumar and J.~R. Partington.
\newblock Weighted composition operators on {H}ardy and {B}ergman spaces.
\newblock In {\em Recent advances in operator theory, operator algebras, and
  their applications}, volume 153 of {\em Oper. Theory Adv. Appl.}, pages
  157--167. Birkh\"{a}user, Basel, 2005.

\bibitem{Le}
T.~Le.
\newblock Normal and isometric weighted composition operators on the {F}ock
  space.
\newblock {\em Bull. Lond. Math. Soc.}, 46(4):847--856, 2014.

\bibitem{littlewood}
J.~E. Littlewood.
\newblock On inequalities in the theory of functions.
\newblock {\em Proc. London Math. Soc. (2)}, 23:481--519, 1925.

\bibitem{LP61}
G.~{Lumer} and R.~S. {Phillips}.
\newblock {Dissipative operators in a Banach space}.
\newblock {\em {Pac. J. Math.}}, 11:679--698, 1961.

\bibitem{mat93}
V.~Matache.
\newblock On the minimal invariant subspaces of the hyperbolic composition
  operator.
\newblock {\em Proc. Amer. Math. Soc.}, 119(3):837--841, 1993.

\bibitem{Mat08}
V.~Matache.
\newblock Weighted composition operators on {$H^2$} and applications.
\newblock {\em Complex Anal. Oper. Theory}, 2(1):169--197, 2008.

\bibitem{mor95}
R.~Mortini.
\newblock Cyclic subspaces and eigenvectors of the hyperbolic composition
  operator.
\newblock In {\em Travaux math\'ematiques, Fasc.\ VII}, S\'em. Math.
  Luxembourg, pages 69--79. Centre Univ. Luxembourg, Luxembourg, 1995.

\bibitem{NRW87}
E.~Nordgren, P.~Rosenthal, and F.~S. Wintrobe.
\newblock Invertible composition operators on {$H^p$}.
\newblock {\em J. Funct. Anal.}, 73(2):324--344, 1987.

\bibitem{pazy}
A.~Pazy.
\newblock {\em Semigroups of linear operators and applications to partial
  differential equations}, volume~44 of {\em Applied Mathematical Sciences}.
\newblock Springer-Verlag, New York, 1983.

\bibitem{pomm}
C.~Pommerenke.
\newblock {\em Boundary behaviour of conformal maps}, volume 299.
\newblock Springer {G}rundlehren, 1993.

\bibitem{rota59}
G.~C. Rota.
\newblock Note on the invariant subspaces of linear operators.
\newblock {\em Rend. Circ. Mat. Palermo (2)}, 8:182--184, 1959.

\bibitem{rota60}
G.-C. Rota.
\newblock On models for linear operators.
\newblock {\em Comm. Pure Appl. Math.}, 13:469--472, 1960.

\bibitem{seyoum}
W.~Seyoum and T.~Mengestie.
\newblock Spectrums and uniform mean ergodicity of weighted composition
  operators on {F}ock spaces.
\newblock {\em Bull. Malays. Math. Sci. Soc.}, 45(1):455--481, 2022.

\bibitem{shap87}
J.~H. {Shapiro}.
\newblock {The essential norm of a composition operator}.
\newblock {\em {Ann. Math. (2)}}, 125:375--404, 1987.

\end{thebibliography}

\end{document}